\documentclass[8pt]{article}
\usepackage{mathrsfs}
\usepackage{bbm}
\usepackage{amsmath}
\usepackage{cases} %
\usepackage{color}

\newtheorem{remark}{Remark}
\newcommand{\bull}{\vrule height 1.8ex width 1.0ex depth 0ex}

\usepackage{graphicx}

\begin{document}{\small
\begin{center}
{\Large\bf On Fast Implementation of Clenshaw-Curtis and Fej\'{e}r-type Quadrature Rules}

\vspace{0.1in} Shuhuang Xiang$^1$, Guo He\footnote{School of Mathematics and Statistics, Central South University, Changsha, Hunan 410083, P. R.
China. Email: xiangsh@mail.csu.edu.cn.  This paper is supported by
National Natural Science Foundation of China No. 11371376.} and Haiyong Wang\footnote{School of Mathematics and Statistics, Huazhong University of Science and Technology, Wuhan, Hubei 430074, P. R. China.}


\end{center}

\vspace{0.2in} {\bf Abstract}. Based upon the fast computation of the coefficients of the interpolation polynomials at Chebyshev-type points by FFT,
DCT and IDST, respectively, together with the efficient evaluation of the modified moments by forwards recursions or by the Oliver's algorithm, this
paper presents interpolating integration algorithms, by using the coefficients and modified moments, for Clenshaw-Curtis, Fej\'er's first and
second-type rules for Jacobi or Jacobi weights multiplied by a logarithmic function. The corresponding {\sc Matlab} codes are included. Numerical
examples illustrate the stability, accuracy of the Clenshaw-Curtis, Fej\'{e}r's first and second rules, and show that the three quadratures have nearly the same convergence rates as Gauss-Jacobi quadrature for functions of finite regularities for Jacobi weights, and are more efficient upon the cpu time than the Gauss evaluated by fast computation of the weights and nodes by {\sc Chebfun}.

 \vspace{0.02in}


\vspace{0.02in}


\vspace{0.25cm}
\section{Introduction}
The interpolation quadrature of  Clenshaw-Curtis rules as well as of the Fej\'{e}r-type formulas for
\begin{equation}
I[f]=\int_{-1}^1f (x)w(x)dx\approx \sum_{k=0}^Nw_k f(x_k):=I_N[f]
\end{equation}
have been extensively studied since Fej\'{e}r \cite{Fejer1,Fejer2} in 1933 and Clenshaw-Curtis \cite{Clenshaw} in 1960, where the nodes $\{x_k\}$ are
of Chebyshev-type while the weights $\{w_k\}$ are computed by sums of trigonometric functions.

\begin{itemize}
\item \textbf{Fej\'{e}r's first-type rule} uses the zeros of the Chebyshev polynomial $T_{N+1}(x)$ of the first kind
$$
y_j=\cos\left(\frac{2j+1}{2N+2}\pi\right),\quad w_j =\frac{1}{N+1}\left\{M_0+2\sum_{m=1}^N M_m\cos\left(m\frac{2j+1}{2N+2}\pi\right)\right\}
$$
for $j=0,1,\ldots,N$, where $\{y_j\}$ is called Chebyshev points of first kind and $M_m=\int_{-1}^1w(x)T_m(x)dx$ (\cite[Sommariva]{Sommariva}).

\item \textbf{Fej\'{e}r's second-type rule} uses the zeros of the Chebyshev polynomial $U_{N+1}(x)$ of the second
kind
$$
x_j=\cos\left(\frac{j+1}{N+2}\pi\right),\quad w_j =\frac{2\sin\left(\frac{j+1}{N+2}\pi\right)}{N+2}\sum_{m=0}^N
\widehat{M}_m\sin\left((m+1)\frac{j+1}{N+2}\pi\right)
$$
for $j=0,1,\ldots,N$, where $\{x_j\}$ is called Chebyshev points of second kind or Filippi points and $\widehat{M}_m=\int_{-1}^1w(x)U_m(x)dx$ (\cite[Sommariva]{Sommariva}).

\item \textbf{Clenshaw-Curtis-type quadrature} is to use the Clenshaw-Curtis points
$$
\overline{x}_j=\cos\left(\frac{j\pi}{N}\right), \quad w_j =\frac{2}{N}\alpha_j\sum_{m=0}^N{'}{'}M_m\cos\left(\frac{jm\pi}{N}\right),\quad
j=0,1,\ldots,N,
$$
where the double prime denotes a sum whose first and last terms are halved, $\alpha_0=\alpha_N=\frac{1}{2}$, and $\alpha_j=1$ for $1\leq j\leq N-1$
(\cite[Sloan and Smith]{Sloan}).
\end{itemize}

In the case $w(x)\equiv 1$, a connection between the Fej\'er and Clenshaw-Curtis quadrature rules
and DFTs was given  by Gentleman \cite{Ge72} in 1972, where the Clenshaw-Curtis rule  is implemented with $N+1$ nodes by means of a discrete
cosine transformation. An independent approach along the same lines, unified
algorithms based on DFTs of order $n$ for generating the weights of the two
Fej\'er rules and of the Clenshaw-Curtis rule, was presented in Waldvogel \cite{Waldvogel} in 2006. A streamlined Matlab code
is given as well in \cite{Waldvogel}. In addition, Clenshaw and Curtis \cite{Clenshaw}, Hara and Smith \cite{Hara}, Trefethen \cite{Trefethen1,Trefethen2}, Xiang and Bornemann in \cite{XiangBornemann}, and Xiang \cite{Xiang2,Xiang3}, etc., showed that
the Gauss, Clenshaw-Curtis and Fej\'er quadrature rules are about equally accurate.

More recently, Sommariva \cite{Sommariva},  following Waldvogel \cite{Waldvogel}, showed that
for general weight function $w$, the weights $\{w_k\}$ corresponding to Clenshaw-Curtis, Fej\'{e}r's first and second-type rules can be computed by IDCT (inverse discrete cosine transform) and DST (discrete sine transform) once the weighted modified moments of Chebyshev polynomials of the first and second kind are available, which generalized the
techniques of \cite{Waldvogel} if the modified moments can be rapidly evaluated.

In this paper, along the way \cite[Trefethen]{Trefethen1}, we consider interpolation approaches for Clenshaw-Curtis rules as well as of the Fej\'{e}r's first and second-type formulas, and present {\sc Matlab} codes for
\begin{equation}
I[f]=\int_{-1}^1f (x)w(x)dx
\end{equation}
for $w(x)=(1-x)^{\alpha}(1+x)^{\beta}$ or $w(x)=(1-x)^{\alpha}(1+x)^{\beta}\ln\left(\frac{1+x}{2}\right)$, which can be efficiently calculated by
FFT, DCT and IDST (inverse DST), respectively: Suppose $Q_N[f](x)=\sum_{j=0}^Na_jT_j(x)$ is the interpolation polynomial at $\{y_j\}$ or
$\{\overline{x}_j\}$, then the coefficients $a_j$ can be efficiently computed by FFT \cite{Ge72,Trefethen1} for Clenshaw-Curtis and by DCT for the
Fej\'er's first-type rule, respectively, and then $I_N[f]=\sum_{j=0}^Na_jM_j(\alpha,\beta)$. So is the interpolation polynomial at $\{x_j\}$ in the
form of $Q_N[f](x)=\sum_{j=0}^Na_jU_j(x)$ by IDST for the Fej\'er's second-type rule with $I_N[f]=\sum_{j=0}^Na_j\widehat{M}_j(\alpha,\beta)$. An
elegant {\sc Matlab} code on the coefficients $a_j$ by FFT for Clenshaw-Curtis points can be found in \cite{Trefethen1}. Furthermore, here the
modified moments $M_j(\alpha,\beta)$ and $\widehat{M}_j(\alpha,\beta)$ can be fast computed by forwards recursions or by Oliver's algorithms with
$O(N)$ operations.

Notice that the fast implementation routine based on the weights $\{w_k\}$ or the coefficient $\{a_k\}$ both will involve in fast computation of the
modified moments. In section 2, we will consider algorithms and present {\sc Matlab} codes on the evaluation of the modified moments. {\sc Matlab}
codes for the three quadratures are presented in section 3, and illustrated  by numerical examples  in section 4.


\section{Computation of the modified moments}
Clenshaw-Curtis-type quadratures are extensively studied in a series of papers by Piessens \cite{Piessens,Piessens1} and Piessens and Branders
\cite{Piessens2,Piessens3,Piessens4}. The modified moment $\int_{-1}^1w(x)T_j(x)dx$ can be efficiently evaluated by recurrence formulae  for Jacobi
weights or Jacobi weights multiplied by $\ln((x+1)/2)$ \cite[Piessens and Branders]{Piessens} in most cases.

\begin{itemize}
\item For $w(x)=(1-x)^{\alpha}(1+x)^{\beta}$:
The recurrence  formula for the evaluation
of the modified moments
\begin{equation}\label{moment1}
  M_k(\alpha,\beta)=\int_{-1}^1w(x)T_k(x)dx,\quad w(x)=(1-x)^{\alpha}(1+x)^{\beta}
\end{equation}
by using Fasenmyer's technique is
\begin{equation}\label{compmoment1}{\small \begin{array}{lll}
  &&(\beta+\alpha+k+2)M_{k+1}(\alpha,\beta)+2(\alpha-\beta)M_k(\alpha,\beta)\\
  &&\quad\quad\quad\quad\quad\quad\quad\quad+ (\beta+\alpha-k+2)M_{k-1}(\alpha,\beta)=0\end{array}}
\end{equation}
with
{\footnotesize$$
M_0(\alpha,\beta)=2^{\beta+\alpha+1}\frac{\Gamma(\alpha+1)\Gamma(\beta+1)}{\Gamma(\beta+\alpha+2)},\quad M_1(\alpha,\beta)=2^{\beta+\alpha+1}\frac{\Gamma(\alpha+1)\Gamma(\beta+1)}{\Gamma(\beta+\alpha+2)}\frac{\beta-\alpha}{\beta+\alpha+2}.
$$}
The forward recursion is numerically stable
\cite[Piessens and Branders]{Piessens}, except in two cases:
\begin{eqnarray}
            \alpha>\beta\quad \mbox{and\quad}\beta&=&-1/2,1/2,3/2,\ldots \\
                          \beta>\alpha\quad \mbox{and\quad}\alpha&=&-1/2,1/2,3/2,\ldots
                             \end{eqnarray}

\item For $w(x)=\ln((x+1)/2)(1-x)^{\alpha}(1+x)^{\beta}$:
For
\begin{equation}\label{moment2}
  G_k(\alpha,\beta)=\int_{-1}^1\ln((x+1)/2)(1-x)^{\alpha}(1+x)^{\beta}T_k(x)dx,
\end{equation}
the recurrence formula \cite{Piessens} is
\begin{equation}\label{compmoment2}{\small\begin{array}{lll}
  &&(\beta+\alpha+k+2)G_{k+1}(\alpha,\beta)+2(\alpha-\beta)G_k(\alpha,\beta)\\
  &&\quad + (\beta+\alpha-k+2)G_{k-1}(\alpha,\beta)=2M_{k}(\alpha,\beta)-M_{k-1}(\alpha,\beta)-M_{k+1}(\alpha,\beta)\end{array}}
\end{equation}
with
{\small$$
G_0(\alpha,\beta)=-2^{\beta+\alpha+1}\Phi(\alpha,\beta+1),\quad G_1(\alpha,\beta)=-2^{\beta+\alpha+1}[2\Phi(\alpha,\beta+2)-\Phi(\alpha,\beta+1)],
$$}
where
$$
\Phi(\alpha,\beta)=B(\alpha+1,\beta)[\Psi(\alpha+\beta+1)-\Psi(\beta)],
$$
$B(x,y)$ is the Beta function and $\Psi(x)$ is the Psi function \cite[Abramowitz and Stegun]{Abram}.
The forward recursion is numerically stable the same as for (2.4) except for (2.5) or (2.6) \cite[Piessens and Branders]{Piessens}.
\end{itemize}

Thus, the modified moments can be fast computed by the forward recursions (2.4) or (2.8) except the cases (2.5) or (2.6) (see Table 1).

\vspace{0.36cm}
\textbf{For the weight $(1-x)^{\alpha}(1+x)^{\beta}$ in the cases of (2.5) or (2.6)}: The accuracy of the forward recursion is catastrophic particularly when
$\alpha-\beta\gg 1$ and $n\gg 1$ (also see Table 2): In case (2.5)  the relative errors $\epsilon_n$ of the computed values $M_n(\alpha,\beta)$
obtained by the forward recursion behave approximately as
$$ \epsilon_n\sim n^{2(\alpha-\beta)}, \quad n\rightarrow \infty$$
and in case (2.6) as
$$ \epsilon_n\sim n^{2(\beta-\alpha)}, \quad n\rightarrow \infty.$$
For this case, we use the Oliver's method \cite{Oliver} with two starting values and one end value to compute the modified moments. Let
\begin{equation}
A_N:=\left(
  \begin{array}{cccc}
    2(\alpha-\beta)  & \alpha+\beta+2+0&                       &  \\
    \alpha+\beta+2-1 & 2(\alpha-\beta) & \alpha+\beta+2+1      &  \\
                     & \ddots          & \ddots                & \ddots \\
                     & \alpha+\beta+2-(N-1) & 2(\alpha-\beta)  & \alpha+\beta+2+(N-1) \\
                     &                      & \alpha+\beta+2-N & 2(\alpha-\beta) \\
  \end{array}
\right),
\end{equation}
\begin{equation}
b_N:=\left(
  \begin{array}{ccccc}
  2^{\alpha+\beta+1}\frac{\Gamma(\alpha+1)\Gamma(\beta+1)}{\Gamma(\alpha+\beta+2)}(\alpha-\beta)&0& \cdots & 0 & -(\alpha+\beta+2+N)M_{N+1}
    \end{array}
\right)^T,
\end{equation}
where ``$\cdot^{T}$'' denotes the transpose, then the
modified moments $M$ can be solved by
\begin{equation}\label{Oliver}
  A_NM=b_N,\quad M=(M_0,M_1,\ldots,M_N)^T,
\end{equation}
where $M_{N+1}$ is computed by hypergeometric function \cite{Piessens} for $N\le 2000$,
\begin{equation}
  M_{N+1}=2^{\alpha+\beta+1}\frac{\Gamma(\alpha+1)\Gamma(\beta+1)}{\Gamma(\alpha+\beta+2)}{_3}F_2([N+1,-N-1,\alpha+1],[1/2,\alpha+\beta+2],1).
\end{equation}
Particularly, if $N>2000$, $M_{N+1}$ is computed by the following asymptotic expression. Taking a change of variables $x=\cos(\theta)$ for
(\ref{moment1}), it yields
\begin{equation}
  M_n(\alpha,\beta)=\int_0^{\pi}\varphi(\theta)\theta^{2\alpha+1}(\pi-\theta)^{2\beta+1}\cos(n\theta) d\theta,\notag
\end{equation}
where
\begin{equation}
 \varphi(\theta)=\bigg(\frac{1-\cos(\theta)}{\theta^2}\bigg)^{\alpha+\frac{1}{2}}\bigg(\frac{1+\cos(\theta)}{(\pi-\theta)^2}\bigg)^{\beta+\frac{1}{2}}, \notag
\end{equation}
then it holds that
\begin{equation}\label{asymptotic1}
M_{n}(\alpha,\beta)= 2^{\beta-\alpha}\sum_{k=0}^{m-1}a_k(\alpha,\beta)h(\alpha+k) + (-1)^n2^{\alpha-\beta}\sum_{k=0}^{m-1}a_k(\beta,\alpha)h(\beta+k)
+ O(n^{-2m})
\end{equation}
by means of the Theorem 3 in \cite[Erd\'{e}lyi]{Eedelyi}, in which
\begin{equation}
  h(\alpha)=\cos\big(\pi(\alpha+1)\big)\Gamma(2\alpha+2)n^{-2\alpha-2},\notag
\end{equation}
$$a_0(\alpha,\beta)=1,\ \ a_1(\alpha,\beta)=-\frac{\alpha}{12}-\frac{\beta}{4}-\frac{1}{6},\ \
 a_2(\alpha,\beta)=\frac{1}{120}+\frac{19\alpha}{1440}+\frac{\alpha^2}{288}+\frac{\alpha\beta}{48}+\frac{\beta}{32}+\frac{\beta^2 }{32}$$
 and $$
  a_3(\alpha,\beta)=-\frac{1}{5040}-\frac{\beta}{960}-\frac{107\alpha}{181440}-\frac{\beta^2}{384}-\frac{\alpha^2}{1920}
  -\frac{\beta^3}{384}-\frac{\alpha^3}{10368}-\frac{7\alpha\beta}{2880}-\frac{\alpha^2\beta}{1152}-\frac{\alpha\beta^2}{384}.
  $$

 The Oliver's algorithm can be fast implemented by applying LU factorization (chasing method) with $O(N)$ operations.

\vspace{0.36cm}
In the case (2.6), by $x=-t$ and $T_n(-x)=\left\{\begin{array}{ll}T_n(x),&\mbox{$n$ even}\\-T_n(x),&\mbox{$n$ odd}\end{array}\right.$,
the computation of the moments can be transferred into the case (2.5).

\vspace{0.36cm}
\textbf{In addition, for the weight $w(x)=\ln((x+1)/2)(1-x)^{\alpha}(1+x)^{\beta}$, in the case (2.5)}: The forward recursion (2.8) is also perfectly
numerically stable (see Table 5) even for $\alpha\gg \beta$. However, in the case (2.6), the forward recursion (2.8) collapses, which can be fixed up
by the Oliver's algorithm similar to (2.9) with two starting values $G_0(\alpha,\beta)$, $G_1(\alpha,\beta)$ and one end value
$G_{N+10^3}(\alpha,\beta)$, by solving an $(N+10^3+1)\times (N+10^3+1)$ linear system for the first $N+1$ moments. The end value can be calculated by
its asymptotic formula, by a change of variables $x=\cos(\theta)$ for (\ref{moment2}) and using
$\ln(\frac{1+\cos(\theta)}{2})=\ln(\frac{1+\cos(\theta)}{2(\pi-\theta)^2})+2\ln(\pi-\theta)$, together with the Theorems in
\cite[Erd\'{e}lyi]{Eedelyi2}, as
\begin{align}\label{asymptotic2}
G_n(\alpha,\beta)= 2^{\beta-\alpha}\sum_{k=0}^{m-1}c_k h(\alpha+k) +
(-1)^n2^{\alpha-\beta}\sum_{k=0}^{m-1}h(\beta+k)\bigg(2a_k(\beta,\alpha)\phi(\beta+k)+b_k\bigg) + O(n^{-2m}),
\end{align}
where
\begin{equation}
\phi(\beta)=\Psi(2\beta+2)-\ln(2n)-\frac{\pi}{2}\tan(\pi\beta),\notag
\end{equation}
and
$$\left\{\begin{array}{l}b_0=0,\ \ b_1=-\frac{1}{12},\ \ b_2= \frac{19}{1440}+\frac{\alpha}{48}+\frac{\beta}{144}, \\
  b_3=-\frac{7\alpha}{2880}-\frac{\beta}{960}-\frac{\alpha^2}{384}-\frac{\beta^2}{3456}-\frac{107}{181440}-\frac{\alpha\beta}{576}, \end{array}\right.
 \left\{ \begin{array}{l}c_0=0,\   c_1=-\frac{1}{4},\ \ c_2= \frac{1}{32}+\frac{\alpha}{48}+\frac{\beta}{16}, \\
  c_3=-\frac{7\alpha}{2880}-\frac{\beta}{192}-\frac{\alpha^2}{1152}-\frac{\beta^2}{128}-\frac{1}{960}-\frac{\alpha\beta}{192}.
\end{array}\right.
$$

\vspace{0.36cm}
Tables 3-6 show the accuracy of the Oliver's algorithm for different $(\alpha,\beta)$, and Table 7 shows the cpu time for implementation of the two
Oliver's algorithms. Here, Oliver-1 means that the Oliver's algorithms with the end value computed by one term of asymptotic expansions, while Oliver-4
signifies that the end value is calculated by four terms of asymptotic expansions. The Oliver-4 can also be applied to the case (2.5) for the Jacobi weight multiplied by $\ln((x+1)/2)$, which can be seen from Table 5 (the Oliver-4 is better than the forward recursion (2.8) in the case (2.5)).

\vspace{0.36cm}
The {\sc Matlab} codes on the Oliver's algorithms and all the {\sc Matlab} codes in this paper can be downloaded from \cite{XiangHeWang}. The all
codes and numerical experiments in this paper are implemented in a Lenovo computer with Intel Core 3.20GHz and 3.47GB Ram.

{\footnotesize\begin{table}[!h] \tabcolsep 0pt \caption{\footnotesize Computation of $M_n(\alpha,\beta)$ and $G_n(\alpha,\beta)$ with different $n$
and $(\alpha,\beta)$ by the forward recursion (2.4) and (2.8) respectively } \vspace*{-15pt}
\begin{center}{\footnotesize
\def\temptablewidth{1.0\textwidth}
{\rule{\temptablewidth}{1pt}}
\begin{tabular*}{\temptablewidth}{@{\extracolsep{\fill}}ccccc}
 n & 10 & 100 & 1000  & 2000   \\    \hline

$\begin{array}{c}\mbox{Exact value for}\\M_n(-0.6,-0.5)\end{array}$ & 0.061104330977316  & 0.009685532923886  & 0.001535055343264  & 0.000881657781753 \\
$\begin{array}{c}\mbox{Approximation by (2.4)}\\\mbox{for $M_n(-0.6,-0.5)$} \end{array}$   & 0.061104330977316  & 0.009685532923886  & 0.001535055343264  & 0.000881657781753 \\
\hline
$\begin{array}{c}\mbox{Exact value for}\\ G_n(10,-0.6)\end{array}$  & -3.053192383855787 & -0.608068551015233  &-0.116362906567503  & -0.070289926350902 \\
$\begin{array}{c}\mbox{Approximation by (2.8) }\\ \mbox{for $G_n(10,-0.6)$}\end{array}$ &-3.053192383855788 & -0.608068551015233 & -0.116362906567506 &  -0.070289926350899

\end{tabular*}
{\rule{\temptablewidth}{1pt}}}
\end{center}
\end{table}}

{\footnotesize\begin{table}[!h]
\tabcolsep 0pt \caption{\footnotesize Computation of $M_n(\alpha,\beta)=\int_{-1}^1(1-x)^{\alpha}(1+x)^{\beta}T_n(x)dx$ with different $n$ and $(\alpha,\beta)$} \vspace*{-15pt}
\begin{center}{\footnotesize
\def\temptablewidth{1\textwidth}
{\rule{\temptablewidth}{1pt}}
\begin{tabular*}{\temptablewidth}{@{\extracolsep{\fill}}cccc}
 n & 5 & 10&100 \\    \hline



\mbox{Exact value for (20,-0.5)} & -1.734810854604316e+05 & 4.049003666168904e+03 & -3.083991348593134e-41 \\

\mbox{(2.4) for (20,-0.5)}     & -1.734810854604308e+05 & 4.049003666169083e+03 & 1.787242305340324e-11 \\    \hline

\mbox{Exact value for (100,-0.5)} & -2.471295049468578e+29 & 1.174275526131223e+29 & 2.805165440968788e-29 \\

\mbox{(2.4) for (100,-0.5)}   & -2.471295049468764e+29 & 1.174275526131312e+29 & -1.380038973213404e+13

\end{tabular*}
{\rule{\temptablewidth}{1pt}}}
\end{center}
\end{table}}

{\footnotesize\begin{table}[!h] \tabcolsep 0pt \caption{\footnotesize Computation of
$M_n(\alpha,\beta)=\int_{-1}^1(1-x)^{\alpha}(1+x)^{\beta}T_n(x)dx$ with $(\alpha,\beta)=(100,-0.5)$ and different $n$ by the Oliver's algorithm}
\vspace*{-15pt}
\begin{center}{\footnotesize
\def\temptablewidth{1\textwidth}
{\rule{\temptablewidth}{1pt}}
\begin{tabular*}{\temptablewidth}{@{\extracolsep{\fill}}ccccc}
 n &2000&4000&8000 \\    \hline

\mbox{Exact value for (0.6,-0.5)}          &  9.551684021848334e-12 & 1.039402748103725e-12   & 1.131065744497495e-13 \\

\mbox{Oliver-4 for (0.6,-0.5)}             &  9.551684021848822e-12 & 1.039402748103918e-12   &  1.131065744497332e-13 \\

\mbox{Oliver-1 for (0.6,-0.5)}             &  9.551684556954339e-12 & 1.039402779428674e-12   &  1.131065757767465e-13 \\ \hline

\mbox{Exact value for (10,-0.5)}          &  -8.412345942129556e-57 & -2.005493070382270e-63   & -4.781368848995069e-70 \\

\mbox{Oliver-4 for (10,-0.5)}             &  -8.412345942129623e-57 & -2.005493070382302e-63   &  -4.781368848995179e-70 \\

\mbox{Oliver-1 for (10,-0.5)}             &  -8.412346024458534e-57 & -2.005493396462483e-63   &  -4.781371046406760e-70

\end{tabular*}
{\rule{\temptablewidth}{1pt}}}
\end{center}
\end{table}}


{\footnotesize\begin{table}[!h] \tabcolsep 0pt \caption{\footnotesize Computation of
$G_n(\alpha,\beta)=\int_{-1}^1(1-x)^{\alpha}(1+x)^{\beta}\ln((1+x)/2)T_n(x)dx$ with different $n$ and $(\alpha,\beta)$ by the Oliver's algorithm}
\vspace*{-20pt}
\begin{center}{\footnotesize
\def\temptablewidth{1.0\textwidth}
{\rule{\temptablewidth}{1pt}}
\begin{tabular*}{\temptablewidth}{@{\extracolsep{\fill}}ccccccc}

n & 10 & 100  & 500   \\    \hline

\mbox{Exact value for (-0.4999,-0.5) }     &-0.314181354550401  &-0.031418104511487&-0.006283620842004 \\

\mbox{Oliver-4 for (-0.4999,-0.5) }   &-0.314181354550428  &-0.031418104511490 & -0.006283620842004 \\

\mbox{Oliver-1 for (-0.4999,-0.5) }   &-0.314181354550438  &-0.031418104511491 & -0.006283620842004 \\ \hline

\mbox{Exact value for (0.9999,-0.5)}       & -0.895286620533541 &-0.088858164406923   & -0.017770353274330 \\

\mbox{Oliver-4  for (0.9999,-0.5) }   & -0.895286620533558 & -0.088858164406925  & -0.017770353274330 \\

\mbox{Oliver-1 for (0.9999,-0.5) }   & -0.895285963133892 &  -0.088858109433133  & -0.017770347359161
\end{tabular*}
{\rule{\temptablewidth}{1pt}}}
\end{center}
\end{table}}

{\footnotesize\begin{table}[!h] \tabcolsep 0pt \caption{\footnotesize Computation of
$G_n(\alpha,\beta)=\int_{-1}^1(1-x)^{\alpha}(1+x)^{\beta}\ln((1+x)/2)T_n(x)dx$ with $(\alpha,\beta)=(100,-0.5)$ and different $n$ by the Oliver's
algorithm} \vspace*{-15pt}
\begin{center}{\footnotesize
\def\temptablewidth{1\textwidth}
{\rule{\temptablewidth}{1pt}}
\begin{tabular*}{\temptablewidth}{@{\extracolsep{\fill}}ccccc}
 n &100&500&1000 \\    \hline

\mbox{Exact value for (100,-0.5)}  & -5.660760361182362e+28 &  -1.126631188200461e+28    &  -5.632306274999927e+27 \\
\mbox{Oliver-4 for (100,-0.5)}     & -5.660760361182364e+28 &  -1.126631188200460e+28    &  -5.632306274999938e+27 \\
\mbox{Oliver-1 for (100,-0.5)}     & -5.660525683370006e+28 &  -1.126606059170211e+28    &  -5.632235588089685e+27 \\
\mbox{(2.8) for (100,-0.5)}        & -5.660760361182770e+28 &  -1.126631188200544e+28    &  -5.632306275000348e+27 \\

\end{tabular*}
{\rule{\temptablewidth}{1pt}}}
\end{center}
\end{table}}

{\footnotesize\begin{table}[!h] \tabcolsep 0pt \caption{\footnotesize Computation of
$G_n(\alpha,\beta)=\int_{-1}^1(1-x)^{\alpha}(1+x)^{\beta}\ln((1+x)/2)T_n(x)dx$ with $(\alpha,\beta)=(-0.5,100)$ and different $n$ by the Oliver's
algorithm compared with that computed by the forward recursion (2.8)} \vspace*{-20pt}
\begin{center}{\footnotesize
\def\temptablewidth{1\textwidth}
{\rule{\temptablewidth}{1pt}}
\begin{tabular*}{\temptablewidth}{@{\extracolsep{\fill}}ccccc}
 n &100&500&1000 \\    \hline

\mbox{Exact value for (-0.5,100)}  &  1.089944378602585e-28   &  7.222157005510106e-198   &  5.715301877322031e-259 \\
\mbox{Oliver-4 for (-0.5,100)}     &  1.089944378602615e-28   &  7.222157005510282e-198   &  5.715301877322160e-259\\
\mbox{Oliver-1 for (-0.5,100)}     &  1.089944378602615e-28   &  7.222157005510282e-198   &  5.715301877322160e-259\\
\mbox{(2.8) for (-0.5,100)}        &  -5.331299059334499e+14  &  -1.061058894110758e+14   &  -5.304494050667818e+13
\end{tabular*}
{\rule{\temptablewidth}{1pt}}}
\end{center}
\end{table}}

{\footnotesize\begin{table}[!h] \tabcolsep 0pt \caption{\footnotesize The cpu time for calculation of the modified moments by the Oliver-4 method for
$\alpha=-0.5$ and $\beta=100$} \vspace*{-10pt}
\begin{center}{\footnotesize
\def\temptablewidth{0.68\textwidth}
{\rule{\temptablewidth}{1pt}}
\begin{tabular*}{\temptablewidth}{@{\extracolsep{\fill}}ccccc}
modified moments&$N=10^3$&$N=10^4$ &$N=10^5$ & $N=10^6$  \\    \hline

$\{M_n(\alpha,\beta)\}_{n=0}^N$&  0.004129s  &    0.012204s  &    0.120747s  &    1.119029s  \\
$\{G_n(\alpha,\beta)\}_{n=0}^N$&  0.006026s  &    0.029010s  &    0.295988s  &    2.902172s  \\


\end{tabular*}
{\rule{\temptablewidth}{1pt}}}
\end{center}
\end{table}}

 \newpage
The {\sc Matlab} codes for weights $M_n(\alpha,\beta)$ and $G_n(\alpha,\beta)$  are as follows:
\begin{itemize}
\item A {\sc Matlab} code for weight $M_n(\alpha,\beta)=\int_{-1}^1(1-x)^{\alpha}(1+x)^{\beta}T_n(x)dx$
$$
\begin{array}{ll}
\texttt{function M=momentsJacobiT(N,alpha,beta)\hspace{1.2cm}\% (N+1) modified moments on $T_n$}&\\
\texttt{f(1)=1;f(2)=(beta-alpha)/(2+beta+alpha);\hspace{1cm}\% initial values}& \\
\texttt{for k=1:N-1}&\\
\texttt{  f(k+2)=1/(beta+alpha+2+k)*(2*(beta-alpha)*f(k+1)-(beta+alpha-k+2)*f(k));}&\\
 \texttt{end;} &  \\
\texttt{M=2\^{}(beta+alpha+1)*gamma(alpha+1)*gamma(beta+1)/gamma(alpha+beta+2)*f;}& \\
\end{array}
$$

\item A {\sc Matlab} code for weight $G_n(\alpha,\beta)=\int_{-1}^1(1-x)^{\alpha}(1+x)^{\beta}\log((1+x)/2)T_n(x)dx$
$$
\begin{array}{ll}
\texttt{function G=momentslogJacobiT(N,alpha,beta)\hspace{.5cm}\% (N+1) modified moments on $T_n$}&\\
\texttt{M=momentsJacobiT(N+1,alpha,beta);\hspace{1cm}\% modified moments on $T_n$ for Jacobi weight}& \\
\texttt{Phi=inline('beta(x+1,y)*(psi(x+y+1)-psi(y))','x','y'); }&\\
\texttt{G(1)=-2\^{}(alpha+beta+1)*Phi(alpha,beta+1);}&\\
\texttt{G(2)=-2\^{}(alpha+beta+2)*Phi(alpha,beta+2)-G(1);}&\\
\texttt{for k=1:N-1}&\\
\texttt{\quad G(k+2)=1/(beta+alpha+2+k)*(2*(beta-alpha)*G(k+1)-}&\\
\texttt{  \quad\quad  \quad\quad \quad        (beta+alpha-k+2)*G(k)+2*M(k+1)-M(k)-M(k+2));}&\\
\texttt{end}&\\
\end{array}
$$
\end{itemize}

The modified moments $\widehat{M}_k(\alpha,\beta)=\int_{-1}^1(1-x)^{\alpha}(1+x)^{\beta}U_k(x)dx$ on Chebyshev polynomials of second kind $U_k$ were considered in Sommariva \cite{Sommariva} by using the formulas
$$
U_n(x)=\left\{\begin{array}{ll}2\sum_{\mbox{$j$ odd}}^nT_j(x),&\mbox{$n$ odd}\\
2\sum_{\mbox{$j$ even}}^nT_j(x)-1,&\mbox{$n$ even}\end{array},\right.
$$
which takes $O(N^2)$ operations for the $N$ moments if $M_k(\alpha,\beta)$ are available. The modified moments $\widehat{M}_k(\alpha,\beta)$ can be efficiently calculated with $O(N)$ operations by using
$$
(1-x^2)U_k'=-kxU_k+(k+1)U_{k-1}
$$
(see Abramowitz and Stegun \cite[pp. 783]{Abram}) and integrating by parts as
\begin{equation}\label{compmoment1}
  (\beta+\alpha+k+2)\widehat{M}_{k+1}(\alpha,\beta)+2(\alpha-\beta)\widehat{M}_k(\alpha,\beta)+ (\beta+\alpha-k)\widehat{M}_{k-1}(\alpha,\beta)=0
\end{equation}
with
$$
\widehat{M}_0(\alpha,\beta)=M_0(\alpha,\beta),\quad \widehat{M}_1(\alpha,\beta)=2M_1(\alpha,\beta),
$$
while for $\widehat{G}_k(\alpha,\beta)=\int_{-1}^1(1-x)^{\alpha}(1+x)^{\beta}\ln((x+1)/2)U_k(x)dx$ as
\begin{equation}\label{compmomentU2}{\small\begin{array}{lll}
  &&(\beta+\alpha+k+2)\widehat{G}_{k+1}(\alpha,\beta)+2(\alpha-\beta)\widehat{G}_k(\alpha,\beta)\\
  &&\quad + (\beta+\alpha-k)\widehat{G}_{k-1}(\alpha,\beta)=2\widehat{M}_{k}(\alpha,\beta)-\widehat{M}_{k-1}(\alpha,\beta)-\widehat{M}_{k+1}(\alpha,\beta)\end{array}}
\end{equation}
with
$$
\widehat{G}_0(\alpha,\beta)=G_0(\alpha,\beta),\quad \widehat{G}_1(\alpha,\beta)=2G_1(\alpha,\beta).
$$

To keep the stability of the algorithms, here we use the following simple equation
\begin{equation}
U_{k+2}=2T_{k+2}+U_k \quad\mbox{(see \cite[pp. 778]{Abram})}
\end{equation}
to derive the modified moments with $O(N)$ operations.
\begin{itemize}

\item A {\sc Matlab} code for weight $\widehat{M}_n(\alpha,\beta)=\int_{-1}^1(1-x)^{\alpha}(1+x)^{\beta}U_n(x)dx$
$$
\begin{array}{ll}
\texttt{function U=momentsJacobiU(N,alpha,beta)\hspace{2.2cm}\% modified moments on $U_n$ }&\\
\texttt{M=momentsJacobiT(N,alpha,beta); \hspace{3.35cm}\% N+1 moments on $T_n$} &\\
\texttt{U(1)=M(1);U(2)=2*M(2);\hspace{5cm}\% initial moments}&\\
\texttt{for k=1:N-1, U(k+2)=2*M(k+2)+U(k); end} &  \\
\end{array}
$$

\item A {\sc Matlab} code for weight $\widehat{G}_n(\alpha,\beta)=\int_{-1}^1(1-x)^{\alpha}(1+x)^{\beta}\log((1+x)/2)U_n(x)dx$
$$
\begin{array}{ll}
\texttt{function U=momentslogJacobiU(N,alpha,beta)\hspace{1.7cm}\% modified moments on $U_n$ }&\\
\texttt{G=momentslogJacobiT(N,alpha,beta);\hspace{3cm}\% modified moments on $T_n$ }& \\
\texttt{U(1)=G(1);U(2)=2*G(2);\hspace{5cm}\% initial moments}\\
\texttt{for k=1:N-1, U(k+2)=2*G(k+2)+U(k); end} &  \\
\end{array}
$$
\end{itemize}


\section{{\sc Matlab} codes for Clenshaw-Curtis and Fej\'{e}r-type quadrature rules}
The coefficients $a_j$ for the interpolation polynomial at $\{\overline{x}_j\}$ can be efficiently computed by FFT \cite{Trefethen1}. For the Clenshaw-Curtis,
we shall not give details but just offer the following {\sc Matlab} functions.

\begin{itemize}
\item For $I[f]=\int_{-1}^1(1-x)^{\alpha}(1+x)^{\beta}f(x)dx$

A {\sc Matlab} code for $I_N^{C\texttt{-}C}[f]$:
$$
\begin{array}{ll}
\texttt{function I=clenshaw{\_}curtis(f,N,alpha,beta)}& \texttt{\% (N+1)-pt C-C quadrature}\\
\texttt{x=cos(pi*(0:N)'/N);} & \texttt{\% C-C points}\\
\texttt{fx=feval(f,x)/(2*N);} & \texttt{\% f evaluated at these points} \\
\texttt{g=fft(fx([1:N+1 N:-1:2]));} & \texttt{\% FFT } \\
\texttt{a=[g(1); g(2:N)+g(2*N:-1:N+2); g(N+1)];} &\texttt{\% Chebyshev coefficients}\\
\texttt{I=momentsJacobiT(N,alpha,beta)*a;} & \texttt{\% the integral}
\end{array}
$$

\item For $I[f]=\int_{-1}^1(1-x)^{\alpha}(1+x)^{\beta}\ln((1+x)/2)f(x)dx$

A {\sc Matlab} code for $I_N^{C\texttt{-}C}[f]$:
 $$
\begin{array}{ll}
\texttt{function I=clenshaw{\_}curtislogJacobi(f,N,alpha,beta)}& \texttt{\% (N+1)-pt C-C quadrature}\\
\texttt{x=cos(pi*(0:N)'/N);} & \texttt{\% C-C points}\\
\texttt{fx=feval(f,x)/(2*N);} & \texttt{\% f evaluated at the points} \\
\texttt{g=fft(fx([1:N+1 N:-1:2]));} & \texttt{\% FFT } \\
\texttt{a=[g(1); g(2:N)+g(2*N:-1:N+2); g(N+1)];} &\texttt{\% Chebyshev coefficients}\\
\texttt{I=momentslogJacobiT(N,alpha,beta)*a;} & \texttt{\% the integral}
\end{array}
$$

\end{itemize}

\textbf{The discrete cosine transform DCT} denoted by $Y={\rm dct}(X)$ is closely related to the discrete Fourier transform but using purely real numbers, and  takes $O(N\log N)$ operations for
{\small$$
Y(k)=w(k)\sum_{s=1}^NX(s)\cos\left(\frac{(k-1)\pi(2s-1)}{2N}\right)\quad \mbox{with $w(1)=\frac{1}{\sqrt{N}}$ and $w(k)=\sqrt{\frac{2}{N}}$ for $2\le k\le N$}.
$$}
\textbf{The discrete sine transform DST} denoted by $Y={\rm dst}(X)$ and its inverse \textbf{The inverse discrete sine transform IDST} denoted  by $X={\rm idst}(Y)$ both takes $O(N\log N)$ operations for
$$
{\small Y(k)=\sum_{s=1}^NX(s)\sin\left(\frac{k\pi s}{N+1}\right).}
$$

Note that the coefficients $a_j$ for the interpolation polynomial ${\displaystyle Q_N(x)=\sum_{j=1}^{N}{'}a_{j-1}T_{j-1}(x)}$ at
$\cos\left(\frac{(2k-1)\pi}{2N}\right)$ are represented by {\small
$$a_{j-1}={\displaystyle \frac{2}{N}\sum_{s=1}^Nf\left(\cos\left(\frac{(2s-1)\pi}{2N}\right)\right)\cos\left(\frac{(2s-1)(j-1)\pi}{2N}\right)},\quad
j=1,2,\ldots,N,
$$}
and $a_j$ for the interpolation polynomial ${\displaystyle Q_N(x)=\sum_{j=1}^{N}a_{j-1}U_{j-1}(x)}$ at $\cos\left(\frac{k\pi}{N+1}\right)$ satisfies
$${\small   f\left(\cos\left(\frac{j\pi}{N+1}\right)\right)\sin\left(\frac{j\pi}{N+1}\right)={\displaystyle
\sum_{s=1}^Na_{s-1}\sin\left(\frac{sj\pi}{N+1}\right)},\quad j=1,2,\ldots,N.}
$$
Then both can be efficiently calculated by DCT and IDST respectively.

\begin{itemize}
\item For $I[f]=\int_{-1}^1(1-x)^{\alpha}(1+x)^{\beta}f(x)dx$

A {\sc Matlab} code for $I_N^{F_1}[f]$:
$$
\begin{array}{ll}
\texttt{function I=fejer1Jacobi(f,N,alpha,beta)}& \texttt{\% (N+1)-pt Fej\'{e}r's first rule}\\
\texttt{x=cos(pi*(2*(0:N)'+1)/(2*N+2));} & \texttt{\% Chebyshev points of 1st kind}\\
\texttt{fx=feval(f,x);} & \texttt{\% f evaluated at these points} \\
\texttt{a=dct(fx)*sqrt(2/(N+1));a(1)=a(1)/sqrt(2);} & \texttt{\% Chebyshev coefficients}\\
\texttt{I=momentsJacobiT(N,alpha,beta)*a;} & \texttt{\% the integral}
\end{array}
$$

A {\sc Matlab} code for $I_N^{F_2}[f]$:
$$
\begin{array}{ll}
\texttt{function I=fejer2Jacobi(f,N,alpha,beta)}& \texttt{\% (N+1)-pt Fej\'{e}r's second rule}\\
\texttt{x=cos(pi*(1:N+1)'/(N+2));} & \texttt{\% Chebyshev points of 2nd kind}\\
\texttt{fx=feval(f,x).*sin(pi*(1:N+1)'/(N+2));} & \texttt{\% f evaluated at these points} \\
\texttt{a=idst(fx);} &\texttt{\% Chebyshev coefficients}\\
\texttt{I=momentsJacobiU(N,alpha,beta)*a;} & \texttt{\% the integral}
\end{array}
$$


\item For $I[f]=\int_{-1}^1(1-x)^{\alpha}(1+x)^{\beta}\ln((1+x)/2)f(x)dx$

A {\sc Matlab} code for $I_N^{F_1}[f]$:
$$
\begin{array}{ll}
\texttt{function I=fejer1logJacobi(f,N,alpha,beta)}& \texttt{\% (N+1)-pt Fej\'{e}r's first rule}\\
\texttt{x=cos(pi*(2*(0:N)'+1)/(2*N+2));} & \texttt{\% Chebyshev points of 1st kind}\\
\texttt{fx=feval(f,x);} & \texttt{\% f evaluated at these points} \\
\texttt{a=dct(fx)*sqrt(2/(N+1));a(1)=a(1)/sqrt(2);} & \texttt{\% Chebyshev coefficients}\\
\texttt{I=momentslogJacobiT(N,alpha,beta)*a;} & \texttt{\% the integral}
\end{array}
$$

A {\sc Matlab} code for $I_N^{F_2}[f]$:
$$
\begin{array}{ll}
\texttt{function I=fejer2logJacobi(f,N,alpha,beta)}& \texttt{\% (N+1)-pt Fej\'{e}r's second rule}\\
\texttt{x=cos(pi*(1:N+1)'/(N+2));} & \texttt{\% Chebyshev points of 2nd kind}\\
\texttt{fx=feval(f,x).*sin(pi*(1:N+1)'/(N+2));} & \texttt{\% f evaluated at these points} \\
\texttt{a=idst(fx);} &\texttt{\% Chebyshev coefficients}\\
\texttt{I=momentslogJacobiU(N,alpha,beta)*a;} & \texttt{\% the integral}
\end{array}
$$

\end{itemize}

\begin{remark}
The coefficients $\{a_j\}_{j=0}^N$ for Clenshaw-Curtis can also be computed by idst, while the coefficients for Fej\'er's rules can be computed by FFT. The following table shows the total time for calculation of the coefficients for $N=10^2:10^4$.

{\footnotesize\begin{table}[!h] \tabcolsep 0pt \caption{\footnotesize Total time  for calculation of the coefficients for $N=10^2:10^4$} \vspace*{-10pt}
\begin{center}{\footnotesize
\def\temptablewidth{0.66\textwidth}
{\rule{\temptablewidth}{1pt}}
\begin{tabular*}{\temptablewidth}{@{\extracolsep{\fill}}ccc}
 \mbox{Clenshaw-Curtis} &\mbox{Fej\'er first} &\mbox{Fej\'er second}  \\    \hline

\mbox{FFT:} 10.539741s &\mbox{FFT:} 16.127888s   &\mbox{FFT:} 9.608675s \\

\mbox{idst:} 12.570079s &\mbox{dct:} 10.449258s   &\mbox{idst:} 10.256482s

\end{tabular*}
{\rule{\temptablewidth}{1pt}}}
\end{center}
\end{table}}
From Table 8, we see that the coefficients computed by the FFT is more efficient than that by the idst for Clenshaw-Curtis, the coefficients computed
by the dct more efficient than that by the FFT for Fej\'er first rule,  and  the coefficients of the interpolant for the second kind of Chebyshev
polynomials $U_n$ computed by the idst nearly equal to that for the first kind of Chebyshev polynomials $T_n$ by the FFT for Fej\'er second rule.
Notice that the FFTs for Fej\'er's rules involves  computation of complex numbers. Here we adopt dct and idst for the two rules.
\end{remark}

\section{Numerical examples}
\setcounter{theorem}{0} \setcounter{equation}{0}
\setcounter{lemma}{0} \setcounter{proposition}{0}
\setcounter{corollary}{0}

The convergence rates of the Clenshaw-Curtis, Fej\'{e}r's first and second rules have been extensively studied in Clenshaw and Curtis
\cite{Clenshaw}, Hara and Smith \cite{Hara}, Riess and Johnson \cite{Riess}, Sloan and Smith \cite{Sloan,Sloan1},  Trefethen
\cite{Trefethen1,Trefethen2}, Xiang and Bornemann in \cite{XiangBornemann}, and Xiang \cite{Xiang2,Xiang3}, etc. In this section, we illustrate the
accuracy and efficiency of the Clenshaw-Curtis, Fej\'{e}r's first and second-type rules for the two functions $\tan|x|$ and $|x-0.5|^{0.6}$ by the
algorithms presented in this paper, comparing with those by the Gauss-Jacobi quadrature used $[x,w]={\rm jacpts}(n,\alpha,\beta)$ in {\sc Chebfun
v4.2} \cite{Trefethen3} (see Figure 1), where the Gauss weights and nodes are fast computed with $O(N)$ operations by Hale and Townsend \cite{Hale} based on Glaser, Liu and V. Rokhlin
\cite{Glaser}. The first column computed by Gauss-Jacobi quadrature in Figure 1 takes 51.959797 seconds and the others
totally take 2.357053 seconds. Additionally, the Gauss-Jacobi quadrature completely fails to compute
$I[f]=\int_{-1}^1(1-x)^{\alpha}(1+x)^{\beta}T_n(x)dx$ for $\alpha\gg 1$ and $n\gg1 $, e.g., $\alpha=100$, $\beta=-0.5$ and $n=100$ (see Table 9). Figure 2 shows
the convergence errors by the three quadrature, which takes  7.336958 seconds.

{\footnotesize\begin{table}[!h] \tabcolsep 0pt \caption{\footnotesize Gauss-Jacobi quadrature $I_n[f]$ for $\int_{-1}^1(1-x)^{100}(1+x)^{-0.5}T_{100}(x)dx$ with $n$ nodes} \vspace*{-10pt}
\begin{center}{\footnotesize
\def\temptablewidth{0.98\textwidth}
{\rule{\temptablewidth}{1pt}}
\begin{tabular*}{\temptablewidth}{@{\extracolsep{\fill}}ccccc}
\mbox{Exact value}& $n=10^2$ &$n=10^3$ &$n=10^4$& $n=10^5$  \\    \hline

2.805165440968788e-29&5.428948613306778e+16&3.412774141453926e+16&8.907453940922673e+17&NaN\\

\end{tabular*}
{\rule{\temptablewidth}{1pt}}}
\end{center}
\end{table}}

\begin{figure}[htbp]
\centerline{\includegraphics[height=9cm,width=16cm]{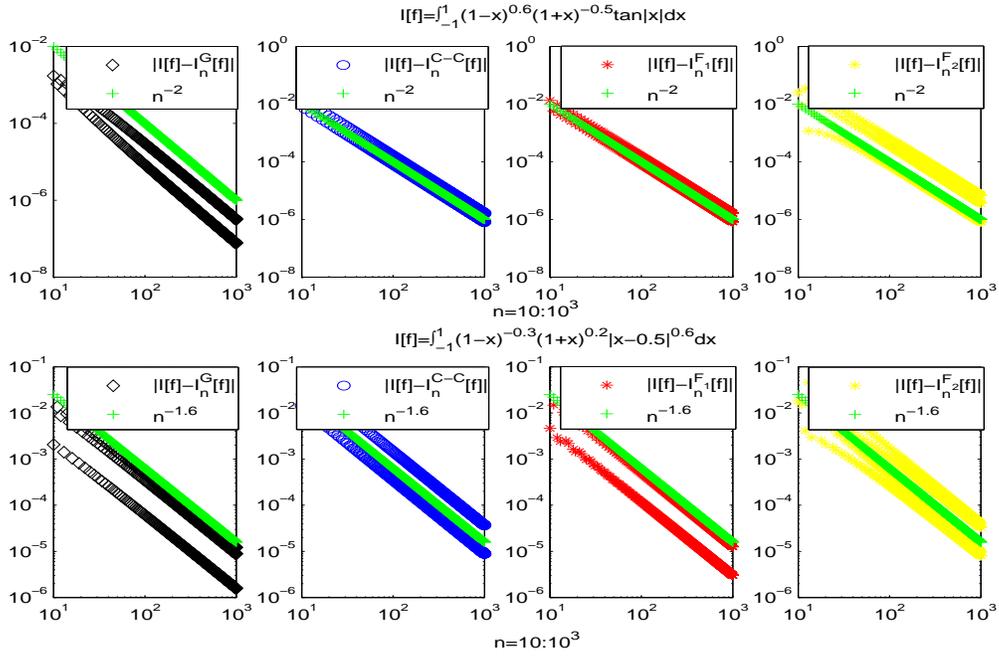}}
\caption{The absolute errors compared with  Gauss quadrature, $n^{-2}$ and $n^{-1.6}$, respectively, for $\int_{-1}^1(1-x)^{\alpha}(1+x)^{\beta}f(x)dx$ evaluated by the Clenshaw-Curtis, Fej\'er's first and second-type rules with $n$ nodes:
$f(x)=\tan|x|$ or $|x-0.5|^{0.6}$ with different $\alpha$ and $\beta$  and $n=10:1000$.}
\end{figure}
\begin{figure}[htbp]
\centerline{\includegraphics[height=9cm,width=16cm]{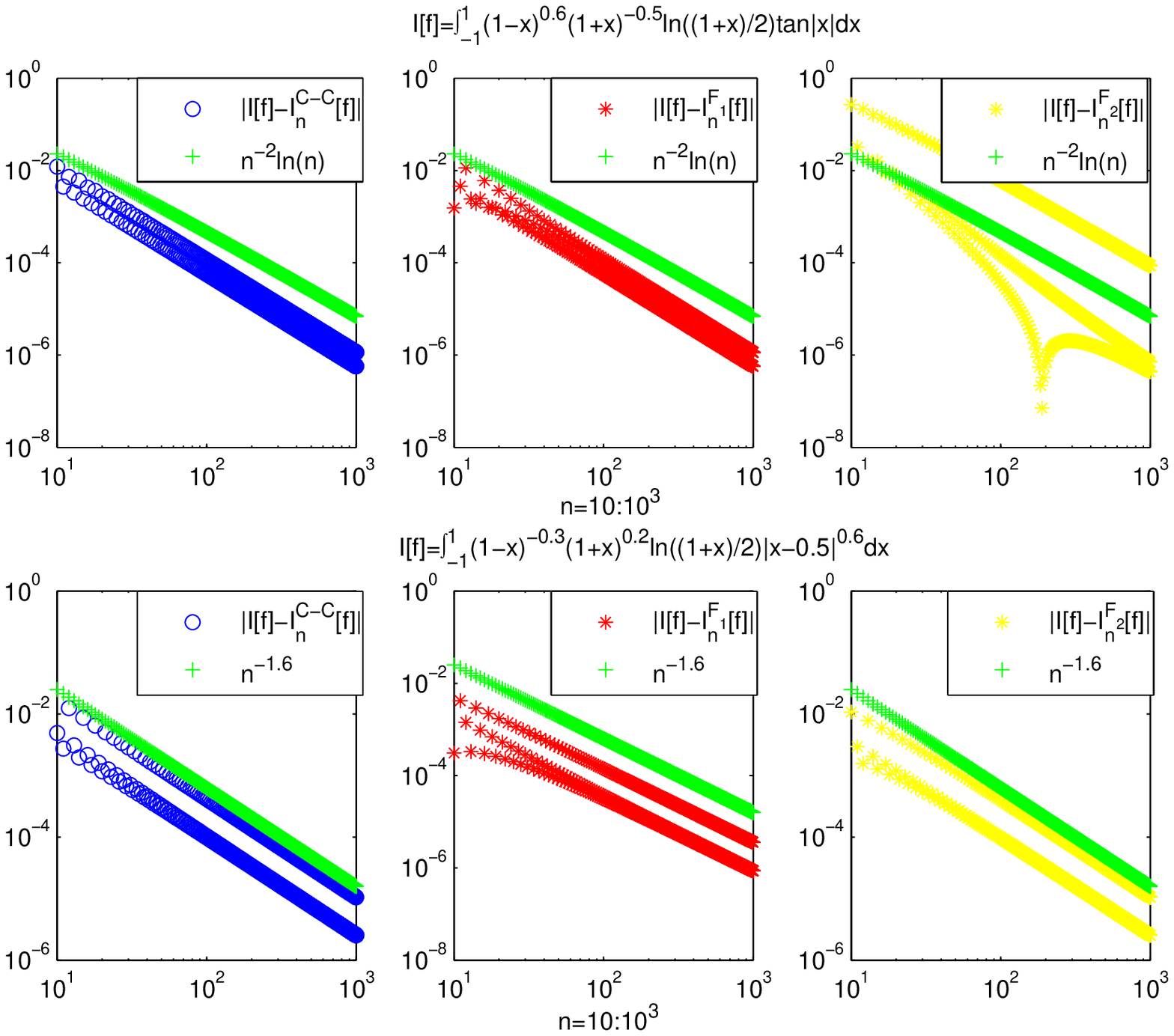}}
\caption{The absolute errors compared with $n^{-2}\ln(n)$ and $n^{-1.6}$, respectively, for $\int_{-1}^1f(x)dx$ evaluated by the Clenshaw-Curtis, Fej\'er's first and second-type reules with $n$ nodes:
$f(x)=\tan|x|$ or $|x-0.5|^{0.6}$ with different $\alpha$ and $\beta$ and $n=10:1000$.}
\end{figure}

Sommariva \cite{Sommariva} showed the efficiency of the computation of the weights $\{w_k\}$ corresponding to Clenshaw-Curtis, Fej\'{e}r's first and second-type rules can be computed by IDCT  and DST for the Gegenbauer weight function
$$w(x) = (1-x)^{\lambda-1/2},\quad\quad \lambda >-1/2
$$
with $\lambda=0.75$ and $N = 2^k$ for $k = 1, \ldots, 20$. Here the modified moments $M_n(\lambda-1/2,\lambda-1/2)$ are available (see \cite[Hunter and Smith]{Hunter}). Table 10 illustrates the cpu time of the computation of the weights $\{w_k\}$ for the computation of
Clenshaw-Curtis, Fej\'er's first and second-type rules by the algorithms given in \cite{Sommariva}, compared with the cpu time of the computation of the coefficients $\{a_k\}$ for the three quadrature by the FFT, DCT and IDST in section 3.

{\footnotesize\begin{table}[!h] \tabcolsep 0pt \caption{\footnotesize The cpu time for calculation of the weight $\{w_k\}_{k=0}^N$ by the algorithms given in \cite{Sommariva}  and the coefficients $\{a_k\}_{k=0}^N$ by the FFT, DCT and IDST in section 3} \vspace*{-10pt}
\begin{center}{\footnotesize
\def\temptablewidth{0.98\textwidth}
{\rule{\temptablewidth}{1pt}}
\begin{tabular*}{\temptablewidth}{@{\extracolsep{\fill}}cccc|cccc}
$\{w_k\}_{k=0}^N$& \mbox{C-C} &\mbox{Fej\'er I} &\mbox{Fej\'er II} & $\{a_k\}_{k=0}^N$& \mbox{C-C} &\mbox{Fej\'er I} &\mbox{Fej\'er II}  \\    \hline

$N=2^{10}$&  0.7152e-3s&    0.4199e-3s&    0.3785e-3s&  $N=2^{10}$&  0.2847e-3s   & 0.3710e-3s   & 0.2905e-3s\\

$N=2^{15}$&  0.0053s&    0.0071s&    0.0087s&  $N=2^{15}$&  0.0052s   & 0.0061s   & 0.0072s\\

$N=2^{18}$&  0.0725s&    0.0871s&    0.1394s&  $N=2^{18}$&  0.0609s   & 0.0604s   & 0.0567s\\

$N=2^{20}$&  0.2170s&    0.2821s&    0.2830s&  $N=2^{20}$&  0.2066s   & 0.2477s   & 0.2345s\\

\end{tabular*}
{\rule{\temptablewidth}{1pt}}}
\end{center}
\end{table}}


\end{document}